\documentclass{amsart}
\title[Prime symmetries]{Surprising symmetries in distribution of prime polynomials}
\author{Dinesh S.\ Thakur}
\address{University of Rochester,  Rochester, NY 14627}
\email{dinesh.thakur@rochester.edu}

\newcommand{\FF}{{\mathbb{F}}}

\newcommand{\X}{\mathcal S}

\begin{document}

\begin{abstract}
The primes or prime polynomials (over finite fields) are supposed to be distributed `irregularly', despite nice asymptotic or average behavior. We provide some conjectures/guesses/hypotheses
with `evidence' of surprising symmetries in prime distribution. 
At least when the characteristic is $2$, we provide 
conjectural rationality and characterization of vanishing for families of interesting infinite sums over irreducible polynomials over finite fields. The cancellations responsible do not happen degree by degree or even for degree bounds for primes or prime powers, so rather than finite fields being 
responsible, interaction between all finite field extensions seems to be playing a role 
 and thus suggests some interesting symmetries in the 
distribution of prime polynomials.  Primes are subtle, so whether there is actual vanishing 
of these sums indicating surprising symmetry (as guessed optimistically), 
or whether these sums just have surprisingly large valuations indicating only 
some small degree cancellation phenomena of low level approximation 
symmetry (as feared sometimes pessimistically!), remains to be seen. In either case, the phenomena begs for an explanation. {\bf A beautiful explanation by David Speyer and important updates are added at the end after the first version.}

\end{abstract}

\maketitle

We start with the simplest case of the conjecture and some of its consequences: 

\vskip .1truein

{\bf Conjecture/Hypothesis  A} When $\wp$ runs through all the primes of $\FF_2[t]$, 

$$\X:= \sum \frac{1}{1+ \wp}=0, \ \ \mbox{as a power series in $1/t$}.$$

\vskip .1truein

The conjecture implies by geometric series development, squaring in characteristic two 
and subtracting 
that 
$$\sum_{n=1}^{\infty} \sum_{\wp} \frac{1}{\wp^n} =0=\sum_{n=1}^{\infty}\sum_{\wp}\frac{1}{\wp^{2n}}
=\sum_{n=1}^{\infty}\sum_{\wp}\frac{1}{\wp^{2n-1}}.$$

Now this series $\X$ is also the logarithmic derivative at $x=1$ of ${\mathcal P}(x) = \prod (1+x/\wp)^{m_{\wp}}$, 
where $m_\wp$'s are any odd integers depending on $\wp$, for example $1$ or $-1$ assigned in arbitrary fashion. 
Hence the derivative of ${\mathcal P}$ at $x=1$ is also zero. 

Now for the simplest choice $m_{\wp}=1$ for all $\wp$, we have 
$${\mathcal P}(x) =\sum_{n=0}^{\infty}\sum \frac{x^n}{\wp_1\cdots\wp_n},$$
where the second sum is over all distinct $n$ primes $\wp_1, \cdots, \wp_n$. 

If we put $m_{\wp}=-1$ for some $\wp$, by the geometric series development of the corresponding term, it contributes to the power series with any multiplicity. 

For such  choices $\pm 1$, this thus implies 
$$0=\sum_{\wp_i} \frac {1}{\wp_1}+\frac{1}{\wp_1\wp_2\wp_3}+\frac{1}{\wp_1\cdots \wp_5}+\cdots,
$$
where the sum is over all the primes with or without multiplicities prescribed for any  subset of primes. 

Choosing $m_{\wp}$'s to be more than one,   bounds the multiplicities of those primes by 
$m_{\wp}$'s, but with complicated sum conditions due to complicated 
vanishing behaviors of the binomial coefficients involved, with $m_{\wp}=2^n-1$ 
giving the simplest nice behavior of just bounding the multiplicity.  

The complementary sum with even number of factors is then (under A) transcendental (zeta value power), at least 
when $m_{\wp}$ is constant. Also note that changing finitely many $m_{\wp}$'s (to even say) 
just changes the answer by a rational function. 

The cancellations happen in a complicated fashion, 
indicating some nice symmetries in the distribution of primes in function fields. Whereas the 
usual sum evaluations involving primes are basically through involvement of all integers through 
the Euler product, here the mechanism seems different, with mysterious Euler product 
connection through logarithmic derivatives. We will explain this later to keep description 
of results and guesses as simple as possible.

Now we explain the set-up and more general conjectures and results.

{\bf Notation}
$$\begin{array}{rcl}

q&=&\mbox{a power of a prime} \ p\\
A&=&\FF_q[t]\\
A+&=&\mbox{\{monics in $A$\}}\\
P & = & \mbox{\{irreducibles in $A$ (of positive degree)\}}\\
P+ & = &  \mbox{\{monic elements in $P$\}}\\
K&=&\FF_q(t)\\
K_\infty&=&\mbox{completion of $K$ at the place $\infty$ of $K$}\\
\mbox{[$n$]}&=&t^{q^n}-t\\
\deg&=&\mbox{function assigning to $a\in A$ its degree in $t$, $\deg(0)=-\infty$}\\
\end{array}
$$

For positive integer $k$, let $P(k) =\sum 1/(1-\wp^k)\in K_\infty$, where the sum is over 
all $\wp\in P+$ and $p(k)$ be its valuation (i.e., minus the degree in $t$) at infinity. 
Similarly we define $P_d(k)$, $p_d(k)$, $P_{\leq d}(k)$, $p_{\leq d}(k)$, by restricting 
the sum to $\wp$'s of degree equal or at most $d$. Since $P(kp)=P(k)^p$, 
we often restrict to $k$ not divisible by $p$. 

\vskip .2truein

{\bf Conjecture / Hypothesis B} For $q=2$, and  $k\geq 1$, 
 $P(k)$ is a rational function in $t$. For example, $P(2^n-1)=[n-1]^2/[1]^{2^n}$.  

We have several  conjectures (the simplest is $k=5$ giving $(t^4+t+1)/([1][3])$) about exact rational functions for hundreds of odd  $k$'s   not of this form, but do not 
have a general conjectural form yet.

\vskip .2truein

{\bf Conjecture/Hypothesis C} Let  $p=2$, and $k$  an odd multiple of $q-1$. Then 
(i)  $P(k)$ is rational 
function in $t$, (ii) When $q=4$, it vanishes if and only if $k=4^n+(4^n-
4^j-1)$, with $n\geq j>0$). (iii) For $q=2^m>2$, it vanishes for $k=2q^n-q^j-1$ with 
$n\geq j >0$  (iv) For $q=2^m$, 
it vanishes, if (and only if ) $p_1(k)\geq 2k$. 

\vfill

\pagebreak

 {\bf Remarks} (0) Again we have several simpler conjectural vanishing criteria for $q=2^m>4$ and 
 conjectural formulas for the rational function, for various $q$'s and $k$'s. {\bf For example},  For $q=4$, $k=9, 21, 57$, it should be  $1/([1]^3+1)$,
$1/[1]^6$ and $[1]/[3]$ respectively, whereas for $q=8, k= 49$, it should be  $[1]/[2]$. 
The `only if' part of (iv) the Conjecture C is trivial (following from trivial weaker form $p_1(k)\geq p_2(k)$).
The `only if' part of  (ii) should be possible to settle, by methods below. But it has not yet been fully worked out.  For $k$'s of part (ii), $P_1(k)=0$, if $n=j$, and we guess $p_1(k) =2k+4^j+2$ otherwise.
 
 (1) It is just possible (but certainly not apparent in our limited computations) 
that rationality works without any restriction on characteristic $p$, but we do not yet have strong evidence either way. Similarly, when $k$ is not a multiple of $q-1$, $P(k)$ does not look 
(we use continued fractions) rational. Is it always transcendental over $K$ then? 

(2) The  valuation results below show that the convergence is approximately linear in the degree
in contrast to usual power sums over all monic polynomials when it is exponential.  This 
makes it hard to compute for large $d$, especially for even moderate size $q$. 

\vskip .1truein

{\bf Sample numerical evidence} Note that the vanishing conjectures and others when the guesses of rational functions are given above, if false,  can be easily refuted by computation. We give only a small sample of computational bounds on accuracy checks that we did by combining calculation 
and theory below (which often improves the bounds a little). 

(A) For $q=2$,  $p(1)\geq 42$ by direct calculation for $d\leq 37$ and some theory. 

(C (ii)) For $q=4$, $p(3)\geq 60, p(15)\geq 228, p(63)\geq 828, p(255)\geq 3072, p(27)\geq 384, p(111)\geq 1224$, with calculation for $d\leq 15, 14, 12, 11, 14, 10$ resp. etc. 

For $q=8$, $p(7)\geq 112, p(21)\geq 224, p(63)\geq 612, p(511)\geq 7168$ etc. 

(B and C(i)): Let $e(k)$ be the valuation of the `error': $P(k)$ minus its guess (so conjecturally infinite). 
For $q=2$, $e(3)\geq 88, e(7)\geq 176, e(15)\geq 348, e(31)\geq 652, e(63)\geq 1324$ etc. 
and $e(5)\geq 130, e(9) \geq 170$ etc. for other guesses. 
For $q=4$, $e(9)\geq 128$ etc. We have a few $q=8, 16$ examples: larger $q$ are hard to compute. 

To illustrate behavior, $q=16$, $k=255$, then $P_1(k)=P_2(k)=0, p_3(k)=3840, p_4(k)=61440, 
p_5(k)=7920$ etc. 

\vskip .1truein
{\bf Guesses/observations at finite levels} Not tried yet to settle. 

(I)  $q=2$,  $p_{2^n}(1)=2^{n+1}+2, n>2$,    $p_{3^n}(1)=3^n+3^{n-1}$,  $p_{5^n}(1)=5^n+5^{n-1}$ ($n\leq 5, 3, 2$ evidence) 
and may be similar for any odd prime power ($n=1$ evidence)?

(II) For general $q$ prime, if $d=2$ or $3$, $p_d(q-1)=q(q-1)$. 
 There is much more data and guesses of this kind, e.g., 
Let $q=4$. If $k=q^{\ell}-1$, then $p_2(k)=p_3(k)=q^{\ell}(q-1)$, 
$p_4(k) = 2q^{\ell}(q-1)$ and $p_{\leq 3}(k)=b_{\ell}$, 
where $b_1=24$ and $b_{n+1}=4b_n+12$. (Evidence $\ell \leq 5$). 

 $q=2^n>2$, $k=q^{\ell}-1$, then $p_3(k)=q^{\ell}(q-1)$ (checked $q=8,16,\ell \leq 4,2$). 

(III) We have guesses for most of (i) $e_{\leq d}(k)$ for $q=2, k=2^n-1, d\leq 17$, e.g. $18k+6$ for $d=17, k>3$. 
(ii) $p_{\leq d}(k)$ for $q=4,  d\leq 10$ where $P(k)$ is guessed zero, e.g., $9k+9$ for $d=7, 8$ and if $k=4^n-1, n>1$, 
also for $d=5$. Also, $q=4, k=4^n-1$, then guess: $p_8(k)=18(k+1), p_4(k)=6(k+1), p_2(k)=3(k+1)$.

\pagebreak

{\bf Remarks}: (1)   Conjecture A is equivalent to : Write $1/(1+\wp)= \sum a_{k, \wp}t^{-k}$, then 
$\sum a_{d, \wp}=0$, where the sum is over all irreducible $\wp$ of degree at most $d$. 
This is equivalent to there being even number of such $\wp$ with $a_{d, \wp}=1$. 

(2)  Note $P(k)$ is  (minus of ) the logarithmic derivative at $x=1$ of (new deformation of Carlitz zeta values)

$$\zeta(x, k): = \sum_{a\in A+} \frac{x^{\Omega(a)}}{a^k}= \prod_{\wp \in P+}(1-\frac{x}{\wp^k})^{-1},$$
where $\Omega(a)$ denotes the number of (monic) prime factors (with multiplicity) of $a$.

The conjecture A was formulated around end 2013-start 2014 via vague optimistic speculations (which could 
not be turned to proofs) about this new zeta variant.  The rationality conjectures were first publicly 
announced at Function Field meeting at Imperial college in summer 2015. All conjectures (except for the explicit forms) follow from another specific conjectural deformation of the Carlitz-Euler formula (for `even' Carlitz zeta values) 
for the zeta variant. Natural candidates seem to have closely related properties.

(3) This is work in progress and in addition to trying to settle these, we are investigating possible generalizations to rationality for general $q$, $L$ functions with characters, logarithmic derivatives at other points, 
higher genus cases and possible number field analogs.  Nothing concrete positive to report yet.
(Except for the following simple observations on easier things at finite level, for whatever  they might be  worth: Since the product of (monic) primes of degree 
dividing $n$ is $[n]$, we get nice formulas for sum of logarithmic derivatives of primes of degree 
$d$, by inclusion-exclusion. We are also looking at power sums $P(d, k)=\sum \wp^k$ 
over degree $d$ monic primes $\wp$.  It seems that the values of degree at most $q$ are only 
$c, c[1]$, with $c\in\FF_p$ (rather than all $c[1]+d$ permitted by translation invariance) (at least in small range of $q, d, k$ checked so far, with 
2 low exceptions for $q=2$), or more generally, if $q>2$ the  product of the constant and $t$-coefficients is zero.
(If true, it should follow from known prescribed coefficients formulae, as we are trying to verify).   A simple sample 
result proved  is $P(d, 1)=0$ for $0<d\leq q-2$, $P(q-1, 1)=-1$, $P(q, 1)=-[1]$ for $q$ prime
(odd for the last part). This follows by simply writing the 
power sum as linear combination (with denominators prime to $q$) of power sums over 
polynomials of degrees at most $d$, which are well-understood. This data has many more patterns
and the first two statements work for $q=4, 8, 9$ also. It is quite possible such things are already somewhere in the literature). 

\vskip .1truein

We hope to put up a more detailed version on ArXives later.

\vskip .2truein

{\bf Acknowledgements} Thanks to John Voight for running a calculation on his excellent computation 
facilities for $q=2, k=1, d\leq 37$ for a couple of months, when I was stuck at $d$ around $20$ using SAGE online. Thanks to my former  student
Alejandro Lara Rodriguez for his help  with SAGE and MAGMA syntax. Thanks to Simon foundation 
initiative which provided MAGMA through UR on my laptop which allowed me to carry out many 
calculations to higher accuracy, much faster 
and often. Thanks to my friends and teachers for their encouragement. 
This research is supported in part 
by NSA grants. 

\vfill

\pagebreak

{\bf Appendix: Sample simple results}: 

\vskip .1truein

{\bf (i) Valuations} For arbitrary $q$ and $k$ a multiple of $q-1$, we have

(0) $P_d(k)\in \FF_p(t^{q-1})$ and it is also a ratio
of polynomials in $[1]$. 

(i) The valuation $p_d(k)$ is divisible by $q(q-1)$ and is at least $kd$. 
(So $p_{\leq d}(k)$ is 
also divisible by $q(q-1)$. )

(ii) We have $p_d(k)= kd$ if and only if (I) $q$ is a prime and $d$ is  square-free multiple of $q$, 
or (II) $q=2$ and $d=4m$ with  $m$ a square-free odd natural number. 

(iii) Let $q=2$ and let $k$ be odd  (without loss of generality) in the usual sense. 
 $P_d(k)$ has $t^{-(dk+1)}$ term if and only if $d$ is square-free. (Hence, 
$p_d(k)=dk+1$ if $d$ is odd square-free,    $p_d(k)=dk$ if $d$ is even square-free 
or $4$ times odd square-free, and $p_d(k)>dk+1$ in other cases). 

\vskip .1truein
The proofs follow by analyzing behavior under automorphisms, and counts of all primes of degree $d$ as well 
as of subset of those containing top two degree terms. 

\vskip .1in

{\bf (ii) Cancellations at finite level}:  When $q=2$, $P_{\leq 2}(1)=0$.
For $q>2$, we have 
$$P_1(q^n-1)=\sum \frac{1}{1-\wp^{q^n-1}}=\sum \frac{\wp}{\wp-\wp^{q^n}}=\frac{\sum (t+\theta)}{t-t^{q^n}}=0. $$

This leads to many more cancellations, unclear whether substantial  for large $d$.

For  $q=2^n>4$, $P_2(q^m-1)=0$. 

\vskip .1truein

This is seen by using above and using that  
when $p=2$, $f(b)=b^2+ab$ is homomorphism from 
$\FF_q$ to itself with kernel $\{0, a\}$

\vskip .1in

{\bf (iii) Non-vanishing of infinite sums}: 
Let $k$ be a multiple of $q-1$ and not divisible by $p$. 

For any $q=2^m$,  if $k>1$  is 1 mod $q$, then $p(k)=p_1(k)=k+(q-1)$, so that $P(k)\neq0$: 

For $q=4$,  $p(k)=p_1(k)$ is the smallest multiple of $12$ greater than $k$ if and only if
$k >3$ and $ k$  is $1$ mod $4$ or $3, 7$ mod $16$. So $P(k)\neq 0$ in these cases. 

\vskip .1truein

These follow by straight analysis of Laurent series expansions using binomial coefficients 
theorem of Lucas.  More results for other $q$'s are in progress. 
\vskip .1truein
\vfill

\pagebreak

{\bf Updates}  After these conjectures were put on polymath blog by Terrence Tao, 
David Speyer gave a very beautiful combinatorial proof for A and of part i of B and C and found the right generalization in any characteristic., using Carlitz' sum and the  product formula for the Carlitz 
exponential and combinatorics of factorization counting and of elementary symmetric functions and power sums. 
See Tao's polymath blog and Speyer's preprint \cite{S}  linked from there. 

Thanks to Terrence Tao and David Speyer! 

\vskip .2truein

We now describe some more results and conjectures based on this progress. 
We use the notation from \cite{S, T}.  

\vskip .2truein

(1) First we remark that using the well-known generalizations of Carlitz exponential 
properties in Drinfeld modules case, Speyer's proof generalizes from $A=\FF_q[t]$ 
case described above to  $A$'s (there are 4 of these, see e.g., \cite [pa. 64, 65]{T}) 
of higher genus with class number one. The general situation, which needs 
a formulation as well as a proof,  is under investigation. 

\vskip .2truein

(2) We have  verified by Speyer's  method a few more isolated conjectural explicit formulas for the rational functions that we had. 
 We have also proved the second part of Conjecture B which gives an explicit family, by 
following Speyer's strategy: 

Let $q=2$. In Speyer's notation, proposition 3.1 of \cite{S}   shows that the 
   the claimed conjecture  is equivalent to
    
    $$g_2(1/a^{2^n-1})=A_n^2, \    \mbox{where}\ \ A_n=\frac{[n-1]}{L_n[1]^{2^{n-1}}}.$$

The left side  is related to power sums by $G_n:=(p_{2^n-1}^2-p_{2(2^n-1)})/2$. 

\vskip .2truein

{\bf Theorem} (I) If we denote by $Y_n$ the reduction modulo $2$ of the standard polynomial expression for $G_n$ in terms of elementary symmetric functions $e_i$ obtained by 
ignoring all monomials which contain $e_i$, with $i$ not of form $2^k-1$, then $Y_n=X_n^2$
with 

$$X_n=\sum_{k=0}^{n-2} e_{2^{n-k}-1}^{2^k} f_k,$$

where $f_0=1$, $f_{k+1}=f_ke_1^{2^k}+X_{k+1}$.

(Equivalently
$$X_n=\sum_{k=0}^{n-2} e_{2^{n-k}-1}^{2^k}(X_k+\sum_{j=1}^{k-1}e_1^{2^{k-1}+\cdots+2^j}X_j)$$
with empty sum being zero convention, the last two terms of bracket could also be combined to get
sum from $0$ to $k$. )

(II) If we substitute $1/D_i$ for $e_{2^i-1}$. and $1/L_i$ for $f_i$ in the formula 
for $X_n$, we get $A_n$. 

(III) In particular, the second part of conjecture B holds. 

\vskip .2truein

{\bf Proof sketch}: 
(I) We have Newton-Girard identities relating power sums to elementary symmetric functions: 

$$p_m=\sum \frac{(-1)^m m (r_1+\cdots r_n-1)!}{r_1!\cdots r_n!}\prod (-e_i)^{r_i}, $$
where the sum is over all non-negative $r_i$ satisfying $\sum ir_i=m$. 

We only care of $i$ of form $N_k:=2^k-1$, as the rest of $e_i$'s are zero, 
when specialized to reciprocals, as in Speyer's proof. Let us put $R_k:= r_{2^k-1}$
and $E_k:=e_{2^k-1}$. Then with $\cong$ denoting ignoring $i$'s not of form $N_k$, we have 

$$p_m\cong \sum \frac{(-1)^m m (R_1+\cdots R_n-1)!}{R_1!\cdots R_n!}\prod (-E_i)^{R_i}, $$
where the sum is over all non-negative $R_k$ satisfying $\sum N_kR_k=m$.

We also only care of $m=2^n-1, 2(2^n-1)$ 
in both of which cases, $k $ can be further be restricted between $1$ and $n$. 

When $m=2^n-1$, we only care about the monomials where the coefficient is 
odd (we do not care about the exact coefficient). This corresponds to the odd multinomial 
coefficients, since  the summing condition  reduced modulo 
2 implies $\sum r_k$ odd, that is why we can reduce to the  multinomial coefficient. 
Now  by Lucas theorem, this corresponds exactly to having 
no clash between the base 2 digits of $R_i$'s.

Multinomial coefficient calculation implies that 
$p_{2^n-1}$ consists of $2^{n-1}$ monomials in $E_k$'s (each of `weight' $ 2^n-1$)
 with odd coefficients, 
out of which (`the second half in lexicographic order') $2^{n-2}$ that make $X_n$  are exactly the ones containing $e_1^s$ with $s\leq 2^{n-1}-1$.  The rest exactly cancel (when squared) with corresponding 
terms from $m=2(2^n-1)$ case (which  we care only modulo $4$) which have odd (or in fact 
3 modulo 4) coefficients, and cross product terms when you square match with the 
even non-zero (modulo 4) coefficients of this case. 
This also explains $Y_n$ is a square and gives formula for its square-root $X_n$.

{\bf Examples}: (The first corresponds to the last example of \cite{S})

$X_2=e_3$

$X_3=e_7       +e_3^2e_1$

$X_4=e_{15} +e_7^2 e_1 +e_3^4(e_3+e_1^3)$

$X_5=e_{31}+e_{15}^2 e_1+e_7^4(e_3+e_1^3)+e_3^8(e_3e_1^4+e_1^7+e_7+e_3^2e_1)$

 $X_6=e_{63}+e_{31}^2 e_1+e_{15}^4(e_3+e_1^3)+e_7^8(e_3e_1^4+e_1^7+e_7+e_3^2e_1)
 +e_3^{16}(....)$

and so on.

(II) (See e.g., \cite{T} for notation, definition, properties) If $\log_C$ and $\exp_C$ denotes the Carlitz logarithm and exponential, which are inverses 
of each other,  $\log_C(\exp_C(z))=z$ gives, by equating coefficients of $z^{2^n}$ for $n>1$, 
 
$$\sum_{k=0}^n 1/(D_{n-k}^{2^k}L_k)=0.$$

 The expression for $X_n$ we get thus reduces to the terms from $k=0$ to $n-2$ whereas 
  $k=n-1, n$ terms give claimed $A_n$, thus proving the claim (as we are in characteristic two). 
  
 The proof is by induction, 
just using 
$$\frac{1}{L_{k+1}}+\frac{1}{L_k[1]^{2^k}}=\frac{1}{L_k}(\frac{1}{[k+1]}+\frac{1}{[1]^{2^{k+1}}})=
\frac{[k]}{L_{k+1}[1]^{2^{k+1}}}$$ 
being the same calculation for $f_{k+1}+f_ke_1^{2^k}$ as well 
as the two terms combination claim above. 

(III) This follows as in \cite{S}, by specializing the symmetric functions to reciprocals of polynomials
and using the evaluation of $e_i$'s coming from Carlitz sum-product formula for Carlitz 
exponential. 

\vskip .2truein

(3) Following \cite{S}, we write 
$ G_p(u)=((1-u^p)-(1-u)^p)/(p(1-u)^p)$. Let us write $G(k)$, $G_d(k)$ etc. for $\sum G_p(1/\wp^k)$. Finite level results similar to (i, ii, iii) on page 5 can be generalized, for example, we have 
valuation divisibility by $q(q-1)$ at finite level, for the same reasons and vanishing for 
degree $1$ and $k=q^n-1$,  for $q$ not a prime. Here is a generalization of the open conjecture 
C (ii, iii) to any $q$. 

\vskip .2truein

(4) {\bf Conjecture / Hypothesis D} Let $k$ be a multiple of $q-1$ and not divisible by $p$. 
(i) Let  $q$ be a prime.  $G(k)=0$ if and only if $k=q-1$, (ii) If $q$ is not a prime, then 
$G(k)=0$ if $k=q^n-1$ or more generally if $k=i(q^n-1)-\sum_{j=1}^r (q^j-1)$, with 
$r<i\leq p$ and $j<n$, (iii) For $q$ a prime, $G(q^n-1)= ([n-1]/[1]^{q^{n-1}})^q$. 

\vskip .2truein

{\bf Remarks} (I) It should be noted that for $q$ of even moderate size (say $9$) the computations blow up quickly, so we do not yet have reasonable/satisfactory  numerical evidence for  D. 
We are still trying to check and improve the statement. 

(II) The `if' part of D (i) is proved in Thm. 1.6 of \cite{S}. The proof of `only if' part 
for $q=2$ in the first version should generalize easily. It is plausible that the techniques and the identities of $q=2$ proof above generalize readily  to settle 
 part (iii) of conjecture D. The combinatorics needed for  part (ii) 
  has not yet been worked out. These  (together with characterization of vanishing) are currently under investigation. 

\vskip .2truein

The author and his  students are working on some of these conjectures, characterization of 
vanishing (which represents extra or clean symmetry),   higher genus and other
variants and the finite level aspects, which will need different tools. 
Till the work is published, more updates will be listed at the end of 

\noindent
http://www.math.rochester.edu/people/faculty/dthakur2/primesymmetryrev.pdf

\end{document}